\documentclass[11pt]{article}
\usepackage[margin=1in]{geometry}
\usepackage{natbib}
\usepackage[symbol]{footmisc} 
\usepackage[hidelinks]{hyperref}
\usepackage{booktabs}
 \usepackage{tabularx}
\usepackage{graphicx}%
\usepackage{eso-pic}
\usepackage{enumitem}
\usepackage{amsfonts,amsmath,amssymb,amsthm,mathtools}%
\usepackage[title]{appendix}%
\usepackage{lmodern}
\usepackage{caption} 

\usepackage{authblk}

\usepackage{times}

\usepackage{tikz}
\usepackage{amsmath}
\usepackage{amssymb}
\usepackage{mathtools}
\usepackage{amsmath}
\usepackage{float}
\usepackage{algorithm}
\usepackage{algpseudocode}
\newtheorem{fact}{Fact}

\renewcommand{\sp}[1]{\|{#1}\|_\text{\rm sp}}
\newcommand{\real}{\mathbb{R}}
\newcommand{\cA}{\mathcal{A}}

\newcommand{\cS}{\mathcal{S}}
\newcommand{\cF}{\mathcal{F}}

\newcommand{\cX}{\mathcal{X}}

\newcommand{\cZ}{\mathcal{Z}}

\renewcommand{\sp}[1]{\|{#1}\|_\text{\rm sp}}

\newcommand{\expec}{\mathbb{E}}

\DeclareMathOperator{\argmin}{argmin}
\newtheorem{assumption}{Assumption}
\newtheorem{theorem}{Theorem}
\newtheorem{remark}{Remark}%
\newtheorem{lemma}{Lemma}
\newtheorem{corollary}[theorem]{Corollary}

\date{}

\begin{document}
\title{\LARGE\bfseries
 Dual-Anchoring Accelerates Stochastic Fixed-Point Iteration 
}


\author{}
\date{}
\maketitle

\vspace{-3em}

\noindent
\begin{minipage}[t]{0.65\textwidth}
    \textbf{Jongmin Lee}\\
    \textit{ University of Pennsylvania}
\end{minipage}%
\hfill
\begin{minipage}[t]{0.33\textwidth}
\raggedleft
    \small  jongmlee@engineering.upenn.edu
\end{minipage}

\begin{abstract}
In this work, we propose the Stochastic Dual-Anchored Krasnoselskii--Mann (SDAKM) iteration for solving stochastic fixed-point problems. SDAKM applies dual-anchoring to stochastic Krasnoselskii--Mann iteration to efficiently reduce the fixed-point residual. In
$(\real^d,\|\cdot\|_p)$ for $p\in[2,\infty]$,
SDAKM exhibits an oracle complexity of
$\widetilde{O}(\epsilon^{-3})$  under iterate-dependent martingale-difference noise with nonexpansive maps. Notably, this improves the state-of-the-art oracle complexity bounds under this noise condition for both Euclidean and non-Euclidean spaces, and  for $p=\infty$, our oracle complexity matches the  lower bound up to logarithmic
factors in $\epsilon$. 
\end{abstract}












    


\section{Introduction}
Given a nonexpansive map $T : \mathbb{R}^d \to \mathbb{R}^d$ with respect to an arbitrary norm $\|\cdot\|$, consider the stochastic fixed-point problem
\[\underset{x\in\real^d}{\text{find }} \,\,\, T(x)=x\]
with access only to stochastic evaluations of $T$. Namely, to find a fixed point of $T$, we can only use evaluations of $T$ corrupted by random noise $U$. In this setup, the stochastic Krasnoselskii-Mann iteration
\begin{equation}
\tag{SKM}
x_{k+1} = (1-\alpha_{k+1})x_{k} + \alpha_{k+1}(Tx_{k} + U_{k+1}), \qquad \forall k\geq 0,
\end{equation}
where $x_0$ is the starting point, $(\alpha_k)_{k\geq 1}$ is a sequence of averaging scalars, and $(U_k)_{k\geq 1}$ is a martingale difference sequence, has been widely studied in the literature on stochastic approximation. Because iterative averaging can reduce the variance of noise through iterations, variants of SKM and their convergence properties have been actively analyzed for solving stochastic fixed-point problems.

However, despite the extensive literature on stochastic approximation, \emph{general} stochastic fixed-point problems have not been thoroughly studied. While variants of SKM achieve near-optimal sample complexity in Hilbert spaces under uniformly bounded noise variance,  oracle complexity in (i) non-Euclidean spaces or (ii) settings with martingale-difference noise without a uniform variance bound has not been fully addressed. A few works study these general settings using SKM, but optimal or near-optimal sample complexity guarantees remain elusive.


Recently, the \emph{Dual Halpern iteration} was introduced as a new sample-efficient accelerated method:
\begin{align*}
x_{k+1}
=x_k+m_k(Tx_k-Tx_{k-1}),  \qquad 0\le k\le K-1,
\end{align*}
 where $x_0$ is the starting point, $Tx_{-1}:=x_0$, and $\{m_k\}_{k=0}^{K-1}$ is a sequence of dual-anchoring coefficients. Here, $m_k(Tx_k-Tx_{k-1})$ is called the \emph{dual-anchoring} term, and in stochastic settings that allow multiple i.i.d. samples, dual-anchoring methods have achieved optimal or near-optimal sample complexity in both Euclidean and non-Euclidean spaces. However, an extension of these results to the martingale-difference setting is not yet known in either Euclidean or non-Euclidean spaces.

\paragraph{Contributions}
In this work, we introduce the Stochastic Dual-Anchored Krasnoselskii--Mann (SDAKM) iteration, which combines SKM and Dual Halpern iteration as follows:
\begin{align*}
x_{k+\frac{1}{2}}
&=(1-\alpha)x_k+\alpha(Tx_k+U_{k+1}),\\
x_{k+1}
&=x_k+m_k(x_{k+\frac{1}{2}}-x_{k-\frac{1}{2}}),
\qquad 0\le k\le K-1,
\end{align*}
 where $x_0$ is the starting point, $x_{-1/2}:=x_0$, and the coefficients $\alpha$ and $\{m_k\}^{K-1}_{k=0}$ depend on the number of iterations $K$.
Under the following iterate-dependent martingale-difference conditions,
\[
\expec[U_{k+1}\mid\cF_k]=0,
\qquad
\expec[\|U_{k+1}\|_p^2\mid\cF_k]
\le a+b\|x_k\|_p^2,
\qquad 0\le k\le K-1,
\]
where $a,b\ge0$, we show that SDAKM achieves the last-iterate rate
\[
\mathbb E\!\left[\|Tx_K-x_K\|_p\right]
=\widetilde{O}(K^{-1/3})
\qquad\text{for}\quad p\in[2,\infty],
\]
where $\widetilde{O}$ hides logarithmic factors in $K$ and other constants. This convergence result directly gives $\widetilde{O}(\epsilon^{-3})$ oracle complexity, and to the best of our knowledge, this improves the state-of-the-art oracle complexity bounds under these noise conditions in both Euclidean spaces ($p=2$)
\cite{tamizholi2026optimal} and
non-Euclidean spaces ($2<p\le\infty$) \cite{bravo2024}. Moreover, for $p=\infty$, our sample complexity is near-optimal in terms of $\epsilon$
by matching the lower bound
in \cite{diakonikolas2026oracle} up to logarithmic factors.

Lastly, we applied SDAKM to stochastic optimization and two-player zero-sum games. Our Stochastic  Dual-Anchored  Gradient Descent and Dual-Anchored $Q$-learning 
achieved $\widetilde{O}(\epsilon^{-3})$ oracle complexity bounds for the  expected gradient norm and the duality gap, respectively. 

\subsection{Prior works}
\paragraph{Stochastic Fixed-point iteration}
The study of stochastic fixed-point iteration has its origins in the context
of statistical estimation and regression \cite{robbins1951stochastic,kiefer1952stochastic}. Depending on the assumptions on the stochastic oracle, various stochastic fixed-point iterations have been developed for solving fixed-point problems \cite{cai2022stochastic, cai2024variance, diakonikolas2026solving, diakonikolas2026oracle}. In stochastic oracle models that allow multiple i.i.d. samples with uniformly bounded variance at chosen query points, \cite{chen2024near} achieves near-optimal $\widetilde{O}(\epsilon^{-2})$ oracle complexity for smooth convex-concave minimax problems in Euclidean spaces. Under the same sampling assumptions, \cite{lee2026optimal} achieves $O(\epsilon^{-3})$ oracle complexity for stochastic fixed-point problems in finite-dimensional non-Euclidean $\ell_p$ spaces with $2<p\le\infty$.
Under martingale-difference noise with uniformly bounded conditional variance in Hilbert spaces, \cite{yang2026stochastic} establishes near-optimal $\widetilde{O}(\epsilon^{-2})$ oracle complexity, matching the lower bound in \cite{foster2019complexity} up to logarithmic factors. Under the same sampling assumptions, \cite{bravo2024} establishes an $O(\epsilon^{-6})$  oracle complexity using stochastic Krasnoselskii--Mann iterations in finite-dimensional normed spaces. Under iterate-dependent martingale-difference noise in finite-dimensional Hilbert spaces, \cite{tamizholi2026optimal} establishes an $O(\epsilon^{-(4+\eta)})$ last-iterate oracle complexity for every fixed $\eta>0$ using stochastic Krasnoselskii--Mann iterations with polynomially decaying stepsizes. For a randomly selected iterate, \cite{cortild2026stochastic} establishes an $O(\epsilon^{-4})$ oracle complexity using stochastic Krasnoselskii--Mann iterations with independent sampling and an additional samplewise nonexpansiveness condition.  In finite-dimensional normed spaces, under iterate-dependent martingale-difference noise and an additional assumption controlling iterate growth, \cite{bravo2024} establishes an $O(\epsilon^{-(10+\eta)})$ last-iterate oracle complexity for every fixed $\eta>0$ using stochastic Krasnoselskii--Mann iterations. 

\paragraph{Krasnoselskii-Mann iteration}
The Krasnosel'ski\u{\i}--Mann iteration was first introduced in the literature on fixed-point problems \cite{mann1953mean, krasnosel1955two} and has been widely studied as a foundational method \cite{bauschke2017correction,moulines2011non, benaim2006dynamics, Watkins1989, littman1994markov}.  In the deterministic setting, its convergence rate has been extensively studied in Hilbert spaces and more general normed spaces \cite{baillon1992optimal, cominetti2014rate}, and in the inexact setting, \cite{bravo2019rates} studied error bounds for iterations of nonexpansive maps. \cite{bravo2024} extended this analysis to the martingale-difference noise setting, and \cite{blaser2026asymptotic} extended it to the Markovian noise setting.

\paragraph{Dual Halpern iteration}
Dual Halpern iteration was first introduced in Hilbert spaces through
the H-duality theorem \cite{yoon2024optimal}, which established optimal complexity in the deterministic setting. In general normed spaces, \cite{lee2026m} also established the optimality of Dual Halpern iteration through the M-duality theorem. In stochastic
oracle models that allow multiple i.i.d. samples with uniformly bounded variance at
chosen query points, \cite{yoon2026direct} established an oracle complexity of
$O(\epsilon^{-3})$ using stochastic Dual Halpern iteration in finite-dimensional Euclidean
spaces under a mean-square nonexpansiveness assumption. \cite{yoon2026dual} further obtained a near-optimal
oracle complexity of $\widetilde{O}(\epsilon^{-2})$ by optimizing the minibatch schedule. Under the same sampling assumptions, \cite{lee2026optimal} established a sample complexity of $O(\epsilon^{-3})$ in finite-dimensional $\ell_p$ spaces. 

\section{Preliminaries}
Throughout the paper, we consider the following assumption on the operator $T$.
\begin{assumption}\label{assum_fixp}
$T:\cX\to\cX$ is a nonexpansive operator on a normed space
$(\cX,\|\cdot\|)$ with a fixed point $x_\star$.
\end{assumption}

For the analysis of the martingale-difference setting, we briefly review the error bound of the inexact Dual Halpern iteration.
\subsection{Error bound of inexact Dual Halpern iteration}
Consider the case in which the operator values $Tx$ can only be approximated by
$\widetilde{T}x=Tx+e$ with approximation error $e$. Define the \emph{Inexact Dual Halpern iteration} \cite{lee2026optimal}  as
\begin{align*}
x^{k+1}    =x^k+m_{K-k}
  (\widetilde{T}x^k-\widetilde{T}x^{k-1}),
\end{align*}
for $k=0,\ldots,K-1$, where
$\widetilde{T}x^k=Tx^k+e_k$,
$\widetilde{T}x^{-1}=x^0$, and $K$ is the total number of iterations. 

\begin{fact}[\cite{lee2026optimal}]\label{fact_idh}
Assume Assumption~\ref{assum_fixp}.
Let $0=m_0\le m_1\le\cdots\le m_K\le1$.
Then the Inexact Dual Halpern iteration exhibits the error bound
\[
\|Tx^K-x^K\|
\le
2(1-m_K+m_KW_K)\kappa
+2\sum_{i=0}^{K-1}W_{K-i}\|e_i\|,
\]
where
\[
W_0=0,
\qquad
W_i=m_i-m_{i-1}+m_{i-1}W_{i-1},
\quad i=1,\ldots,K.
\]
Specifically, if $m_k=k/(k+2)$, then
\[
\|Tx^K-x^K\|
\le
\frac{8\kappa}{K+4}
+4\sum_{i=0}^{K-1}\frac{\|e_i\|}{K-i+2}.
\]
\end{fact}


This fact separately characterizes the convergence rate in terms of $\kappa$ and the accumulated approximation error of the inexact Dual Halpern iteration. Specifically, when $m_k=k/(k+2)$ and $\|e_i\|\le e$, the first term decreases at a rate of $1/K$, and the second term increases at a rate of $\log(K+1)$ as the number of iterations $K$ increases. Inspired by these results, we combine the Dual Halpern iteration with the Krasnosel'skii–Mann iteration to control the stochastic error term through an averaging step while reducing the $\kappa$ term through dual anchoring.

\subsection{Notation}
Write $(\cX,\|\cdot\|)$ for a normed space. An operator $T: \cX \rightarrow \cX$ is nonexpansive if $\|Tx-Ty\| \le \|x-y\|$ for all $x,y \in \cX$.
We denote the $\ell_p$ norms in $\mathbb{R}^d$ by $\|\cdot\|_p$ for $p\in[1,\infty]$, i.e.,
for $x = (x_1,\ldots,x_d) \in \mathbb{R}^d$:
$\|x\|_p = \left(\sum_{i=1}^d |x_i|^p\right)^{1/p}$
for $p \in [1,\infty)$, and
$\|x\|_\infty = \max_{1\le i\le d}|x_i|$.
We say $x_\star$ is a fixed point of the operator $T$ if $Tx_\star=x_\star$, and write
$\kappa=\|x_0-x_\star\|$ for a given $(\cX,\|\cdot\|)$.

\section{Stochastic Dual-Anchored Krasnoselskii-Mann iteration}
We now introduce our novel method, the \emph{Stochastic Dual-Anchored Krasnoselskii-Mann iteration (SDAKM)}:
\begin{align*}
x_{k+\frac{1}{2}}&=(1-\alpha)x_k+\alpha(Tx_k+U_{k+1})\\
x_{k+1}
&=x_k+\frac{K-k}{K-k+2}(x_{k+\frac{1}{2}}-x_{k-\frac{1}{2}})
\qquad 0\le k\le K-1,
\end{align*}
where $K\ge1$ is the number of iterations, $x_0$ is the starting point, $x_{-1/2}:=x_0$, and $U_{k+1}$ is the stochastic error at iteration $k$. Motivated by \cite{lee2026optimal}, we set $m_k=(K-k)/(K-k+2)$, and
$\alpha\in(0,1]$ is a constant stepsize that depends on $K$ and is specified in Theorem~\ref{thm_resi}.

Our SDAKM directly combines SKM and Dual Halpern iteration. 
As we noted, the rough intuition behind SDAKM is that dual-anchoring reduces the term related to $\kappa=\|x_0-x_\star\|$ and averaging reduces the stochastic error related to $U_k$. We will validate this intuition through theoretical analysis in the next subsection.

Throughout the paper, we consider iterate-dependent martingale differences for the stochastic error.

\begin{assumption}\label{assum_stoo}
$(U_k)_{k\ge1}$ is a martingale difference sequence 
with respect to a filtration $(\cF_k)_{k\ge0}$ and satisfies
\[
\expec[U_{k+1} \mid \cF_{k}]=0,
\qquad
\expec[\|U_{k+1}\|^2 \,|\,  \cF_{k}] \le a+b \|x_{k}\|^2
\quad \text{for} \quad 0\le a,b, \quad  0\le k\le K-1.
\]
\end{assumption}
This type of stochastic error, which may depend on the iterate $x_n$ in the approximation of $Tx_n$, appears widely in stochastic optimization, game theory, and reinforcement learning. Nevertheless, stochastic oracle complexity under Assumption~\ref{assum_stoo} has not been thoroughly studied as we clarified in the discussion of prior work.

\subsection{Error bound with averaged martingale differences}

For the analysis of the convergence of SDAKM, we first point out that using a bound for each $U_n$, as in the inexact setup in the Preliminaries section, does not give an efficient sample complexity bound. To use the martingale difference property, specifically Fact~\ref{fact_sam} in the next subsection, we should consider a weighted average of the $U_n$. To obtain such a term, following the strategy of \cite[Section 2.1]{bravo2024}, define $\overline{U}_k$ by
\begin{align*}
\overline{U}_{k+\frac{1}{2}}&=(1-\alpha)\overline{U}_k+\alpha U_{k+1}\\
\overline{U}_{k+1}&=\overline{U}_k+\frac{K-k}{K-k+2}
(\overline{U}_{k+\frac{1}{2}}-\overline{U}_{k-\frac{1}{2}})
\qquad 0\le k\le K-1,
\end{align*}
where $\overline{U}_0=\overline{U}_{-\frac{1}{2}}=0$. 
The advantage of this alternative form of the iteration is that, in contrast with
the original noise $U_n$, the averages $\overline{U}_n$ may converge to zero fast enough. 

The following lemma specifies the averaging coefficients of $\overline{U}_k$.



\begin{lemma}\label{lem_avg_mar} Assume Assumption~\ref{assum_stoo}.
For $1\le k\le K$ and $0<\alpha\le1$, SDAKM gives $\overline{U}_k=\sum^k_{i=1}c_{k,i}U_i$ with deterministic coefficients satisfying
$0<c_{k,i}\le\alpha$ for $1\le i\le k$ and $\sum^k_{j=1}c_{k,j}\le1$.
\end{lemma}
\begin{proof}
By the previously defined iteration for $\overline{U}_k$, we have, for $0\le k\le K-1$,
\[
\overline{U}_{k+1}=m_{k}\overline{U}_{k+\frac{1}{2}}
+\sum^{k-1}_{i=0}(m_{i}-m_{i+1})\overline{U}_{i+\frac{1}{2}},
\]
where $m_i=(K-i)/(K-i+2)$. The weights on the right-hand side are positive and sum to $m_0\le1$.
Then, since $\overline{U}_{k+\frac{1}{2}}=(1-\alpha)\overline{U}_k+\alpha U_{k+1}$,
$0<c_{k,i}\le\alpha$ and $\sum^k_{i=1}c_{k,i}\le1$ follow directly by induction.

\end{proof}

Based on Lemma~\ref{lem_avg_mar}, we can obtain a version of the error bound in Fact~\ref{fact_idh} in terms of the averaged martingale differences $\overline{U}_k$.
\begin{lemma}\label{lem_mt_bd}
Assume Assumptions~\ref{assum_fixp} and~\ref{assum_stoo}. SDAKM exhibits the error bound
\[
\|Tx_K-x_K\|\le\frac{8\kappa}{\alpha(K+4)}
+4\sum_{k=0}^{K-1}\frac{\|\overline{U}_k\|}{K-k+2}+2\|\overline{U}_K\|,
\]
and 
\[
\|x_t-x_\star\|
\le\kappa+\|\overline{U}_t\|
+\alpha\sum_{j=0}^{t-1}\|\overline{U}_{j}\|
\qquad 0\le t\le K.
\]
\end{lemma}
\begin{proof}  Subtracting the iteration for $\overline{U}_k$ from the SDAKM iteration gives
\begin{align*}
x_{k+\frac{1}{2}}-\overline{U}_{k+\frac{1}{2}}
&=(1-\alpha)(x_k-\overline{U}_k)+\alpha Tx_k\\
x_{k+1}-\overline{U}_{k+1}
&=x_k-\overline{U}_k+m_k
\bigl(x_{k+\frac{1}{2}}-\overline{U}_{k+\frac{1}{2}}
-(x_{k-\frac{1}{2}}-\overline{U}_{k-\frac{1}{2}})\bigr)
\qquad 0\le k\le K-1.
\end{align*}
Let $z_k=x_k-\overline{U}_k$ and
$z_{k-\frac{1}{2}}=x_{k-\frac{1}{2}}-\overline{U}_{k-\frac{1}{2}}$.
This gives 
\begin{align*}
z_{k+\frac{1}{2}}&=(1-\alpha)z_k+\alpha T(z_k+\overline{U}_k)\\
z_{k+1}
&=z_k+m_k(z_{k+\frac{1}{2}}-z_{k-\frac{1}{2}})
\qquad 0\le k\le K-1,
\end{align*}
where $z_{-\frac{1}{2}}=z_0=x_0$.
Next, define the operator
$S=(1-\alpha)I+\alpha T$, which is nonexpansive with a fixed point $x_\star$, and let $e_k=\alpha(T(z_k+\overline{U}_k)-Tz_k)$.
Then we get
\begin{align*}
z_1&=z_0+m_0(Sz_0+e_0-z_0),\\
z_{k+1}
&=z_k+m_k(Sz_k+e_k-Sz_{k-1}-e_{k-1})
\qquad 1\le k\le K-1.
\end{align*}
By Fact~\ref{fact_idh} and $\|e_k\|\le\alpha\|\overline{U}_k\|$, we get
\[
\|Sz_K-z_K\|\le\frac{8\kappa}{K+4}
+4\sum_{k=0}^{K-1}\frac{\alpha\|\overline{U}_k\|}{K-k+2}.
\]
Since $Sz-z=\alpha(Tz-z)$, we have
\begin{align*}
\|(I-T)x_k\|
&\le\|(I-T)z_k\|+\|(I-T)x_k-(I-T)z_k\|\\
&\le\|Tz_k-z_k\|+2\|x_k-z_k\|\\
&\le\alpha^{-1}\|Sz_k-z_k\|+2\|\overline{U}_k\|,
\end{align*}
where the second inequality follows from the fact that $I-T$ is $2$-Lipschitz. 
Taking $k=K$ proves the first claim.

We now prove the second claim. First, we have
\[
z_{t+1}
=(1-m_0)x_0+m_t(Sz_t+e_t)
+\sum_{j=0}^{t-1}(m_j-m_{j+1})(Sz_j+e_j),
\]
and $
\|Sz_j+e_j-x_\star\|
\le\|z_j-x_\star\|+\|e_j\|$ by the nonexpansiveness of $S$. Then, by induction on $z_t$, we have
\begin{align*}
\|z_{t+1}-x_\star\|
&\le(1-m_0)\kappa
+m_t\left(\kappa+\sum_{i=0}^{t}\|e_i\|\right)
+\sum_{j=0}^{t-1}(m_j-m_{j+1})
\left(\kappa+\sum_{i=0}^{j}\|e_i\|\right)\\
&\le(1-m_0)\kappa
+\left(m_t+\sum_{j=0}^{t-1}(m_j-m_{j+1})\right)
\left(\kappa+\sum_{i=0}^{t}\|e_i\|\right)\\
&\le\kappa+\sum_{i=0}^{t}\|e_i\|.
\end{align*}
Finally, since $x_t=z_t+\overline{U}_t$ and $\|e_j\|\le\alpha\|\overline{U}_{j}\|$, we get
\[
\|x_t-x_\star\|
\le\|z_t-x_\star\|+\|\overline{U}_t\|
\le\kappa+\|\overline{U}_t\|
+\alpha\sum_{j=0}^{t-1}\|\overline{U}_{j}\|.
\]
\end{proof}

Through this error bound, we now analyze convergence of SDAKM  with respect to $\|\cdot\|_p$.

\subsection{Expected Residual with the \texorpdfstring{$\|\cdot\|_p$ norm}{p-norm}}
In this subsection, we let $(\cX,\|\cdot\|)=(\real^d,\|\cdot\|_p)$
with $p\in(1,\infty]$. The following fact bounds the $r_p$th moment of the norm of a sum
of martingale differences in $(\real^d,\|\cdot\|_p)$.  In particular, for $1<p \le 2$, refer to \cite[Theorem 1.1]{pinelis2015best} with $f(t)=|t|^p$, and for $2<p \le \infty$, refer to \cite[Section 2, Example 3.2]{juditsky2008large}.

\begin{fact}[\cite{pinelis2015best, juditsky2008large}]\label{fact_sam}
Let $p\in(1,\infty]$, and let
$U_1,\ldots,U_m\in\real^d$ be a martingale difference sequence
with respect to a filtration $(\mathcal F_i)_{i=0}^m$, satisfying
\[
\expec[U_i\mid\mathcal F_{i-1}]=0,
\qquad
\expec[\|U_i\|_p^2]<\infty,
\qquad i=1,\ldots,m.
\]
Then, for any deterministic coefficients $a_1,\ldots,a_m\in\real$,
\[
\left(
\expec\left[
\left\|\sum_{i=1}^m a_i U_i\right\|_p^{r_p}
\right]
\right)^{1/r_p}
\le
C_p
\left(
\sum_{i=1}^m
|a_i|^{r_p}
\expec\left[\|U_i\|_p^{r_p}\right]
\right)^{1/r_p},
\]
where 
\[
r_p=
\begin{cases}
\displaystyle p, & 1<p\le2,\\[6pt]
2, & 2<p\le\infty,
\end{cases}
\qquad
C_p=
\begin{cases}
2^{1/p}, & 1<p<2,\\[4pt]
\sqrt{p-1}, & 2\le p<\infty,\\[4pt]
\sqrt{2e\ln(2d)}, & p=\infty.
\end{cases}
\]
\end{fact}

We briefly note that this fact and the results below can be extended to $\kappa$-regular spaces \cite{juditsky2008large}, which include $\ell_p$ spaces.

Combining this fact with the previous lemmas, we now present our main theorem. 

\begin{theorem}\label{thm_resi}
Assume Assumptions~\ref{assum_fixp} and~\ref{assum_stoo}.
For $ 1< p \le \infty $, let
$\alpha=D_p^{-\frac{r_p}{r_p-1}}(KL_K)^{-\frac{r_p}{2r_p-1}}$, where $L_K=2+4\sum_{k=0}^{K-1}(K-k+2)^{-1}$ and
$D_p=\max\{1,3C_p\sqrt{b}\}$. Then, SDAKM exhibits the rate
\[\mathbb E\!\left[\|Tx_K-x_K\|_p\right]
\le
\begin{cases}
\displaystyle
K^{-\frac{p-1}{2p-1}}
\bigl(\ln(K+2)\bigr)^{\frac{p}{2p-1}} Q_p,
&1<p< 2,
\\[12pt]
\displaystyle
K^{-\frac{1}{3}}\bigl(\ln(K+2)\bigr)^{\frac{2}{3}}Q_p,
&2 \le  p \le \infty,
\end{cases}\]
where $Q_p=4^{\frac{r_p}{2r_p-1}}\left(8\kappa D_p^{\frac{r_p}{r_p-1}}+\frac{2C_p}{D_p}\sqrt{a+b(\|x_\star\|_p+\kappa)^2}\right).$
\end{theorem}
\begin{proof}
For readability, write $r=r_p$, $D=D_p$. 
Recall that $\overline{U}_0=0$ and that, by Lemma~\ref{lem_mt_bd},
\[
\|x_k\|_p
\le \|x_\star\|_p+\kappa+\|\overline{U}_k\|_p
+\alpha\sum_{j=0}^{k-1}\|\overline{U}_{j}\|_p.
\]
Since $r\le2$, conditional Jensen's inequality gives
\begin{align*}
\left(
\expec[\|U_{k+1}\|_p^r\mid\mathcal F_k]
\right)^{1/r}
&\le
\left(
\expec[\|U_{k+1}\|_p^2\mid\mathcal F_k]
\right)^{1/2}\\
&\le \sqrt{a+b\|x_k\|_p^2}\\
&\le \sqrt{a+b(\|x_\star\|_p+\kappa)^2}+\sqrt{b}\left(
\|\overline{U}_k\|_p+\alpha\sum_{j=0}^{k-1}\|\overline{U}_{j}\|_p
\right),
\end{align*}
where the second inequality is from Assumption~\ref{assum_stoo} and the last inequality follows from Lemma~\ref{lem_mt_bd} and $
\sqrt{a+b(A+B)^2}\le\sqrt{a+bA^2}+\sqrt{b}\,B$ for $A,B\ge0$. 

Let
\[
C=\sqrt{a+b(\|x_\star\|_p+\kappa)^2}, \qquad \hat{U}_k=\left(\expec[\|\overline{U}_k\|_p^r]\right)^{1/r},
\qquad
U=\max_{0\le k\le K}\hat{U}_k.
\]

Taking $L_r$ norms and using Minkowski's inequality, we obtain
\begin{align*}
\left(\expec[\|U_{k+1}\|_p^r]\right)^{1/r}
&\le C+\sqrt{b}\left(
\hat{U}_k+\alpha\sum_{j=0}^{k-1}\hat{U}_{j}
\right)\\
&\le C+\sqrt{b}(1+\alpha K)U.
\end{align*}

Also, Lemma~\ref{lem_avg_mar} gives
\[
\sum_{i=1}^k c_{k,i}^r
\le \alpha^{r-1}\sum_{i=1}^k c_{k,i}
\le \alpha^{r-1},
\]
and since $\overline{U}_k=\sum_{i=1}^k c_{k,i}U_i$, we have
\begin{align*}
\hat{U}_k
&\le C_p\left(
\sum_{i=1}^k c_{k,i}^r\expec[\|U_i\|_p^r]
\right)^{1/r}\\
&\le C_p\alpha^{(r-1)/r}
\left(C+\sqrt{b}(1+\alpha K)U\right),
\end{align*}
where the first inequality is from Fact~\ref{fact_sam} and the second inequality is from the previous inequalities.
Taking the maximum over $k$ yields
\[
U\le C_pC\alpha^{(r-1)/r}+C_p\alpha^{(r-1)/r}
\sqrt{b}(1+\alpha K) U.
\]
By the definition of $\alpha$ in Theorem~\ref{thm_resi},
\begin{align*}
C_p\alpha^{(r-1)/r}
\sqrt{b}(1+\alpha K)&=\frac{C_p\sqrt{b}}{D}
\left(
(KL_K)^{-\frac{r-1}{2r-1}}+\frac{D^{-\frac{r}{r-1}}}{L_K}
\right)
\\&\le\frac13\left(1+\frac12\right)
=\frac12,
\end{align*}
where the last inequality follows from the definition of $D$, $K\ge1$, and $L_K\ge2$.  Consequently,
\[
U\le2C_pC\alpha^{(r-1)/r}
=\frac{2C_pC}{D}(KL_K)^{-\frac{r-1}{2r-1}}.
\]
 Also, Jensen's inequality gives
$\expec[\|\overline{U}_k\|_p]\le \hat{U}_k\le U$.
Finally, taking expectations
in the residual bound from Lemma~\ref{lem_mt_bd}, we obtain
\begin{align*}
\expec[\|Tx_K-x_K\|_p]
&\le
\frac{8\kappa}{\alpha(K+4)}
+4\sum_{k=0}^{K-1}
\frac{\expec[\|\overline{U}_k\|_p]}{K-k+2}
+2\expec[\|\overline{U}_K\|_p]\\
&\le \frac{8\kappa}{\alpha K}+L_KU\\
&\le
\left(8\kappa D^{\frac{r}{r-1}}+\frac{2C_pC}{D}\right)
K^{-\frac{r-1}{2r-1}}L_K^{\frac{r}{2r-1}}.
\end{align*}

Since $L_K=2+4\sum_{j=3}^{K+2}j^{-1}\le4\ln(K+2)$, the claimed bound follows.
\end{proof}
\begin{corollary}\label{cor_sam}
Assume Assumptions~\ref{assum_fixp} and~\ref{assum_stoo}.
Let $1<p\le\infty$, and let $Q_p$ and $\alpha$ be as in Theorem~\ref{thm_resi}. Then, for $0<\epsilon\le Q_p$, SDAKM gives
\[
\expec[\|Tx_K-x_K\|_p]\le\epsilon
\]
with the number of iterations
\[K=
\begin{cases}
\displaystyle
\widetilde O\left(
\left(\frac{Q_p}{\epsilon}\right)^{\frac{2p-1}{p-1}}
\right),
&1<p < 2,
\\[12pt]
\displaystyle
\widetilde O\left(
\left(\frac{Q_p}{\epsilon}\right)^3
\right),
&2 \le p \le \infty,
\end{cases}
\]
where $\widetilde O$ suppresses logarithmic factors.
\end{corollary}

\begin{proof}
Write $r=r_p$.
By Theorem~\ref{thm_resi},
\[
\expec[\|Tx_K-x_K\|_p]
\le
Q_p\left(
\frac{(\log(K+2))^{r/(r-1)}}{K}
\right)^{\frac{r-1}{2r-1}}.
\]
Choose $M_r\ge1$ and $K$ to satisfy
\[
M_r\ge
\left(
\frac{3r-1}{r-1}+\frac{\log(4M_r)}{\log3}
\right)^{\frac{r}{r-1}}, \qquad K=\left\lceil
M_r\left(\frac{Q_p}{\epsilon}\right)^{\frac{2r-1}{r-1}} \left(\log\left(2+\frac{Q_p}{\epsilon}\right)\right)^{\frac{r}{r-1}}
\right\rceil.
\]
Let $R=\frac{Q_p}{\epsilon}$ and $h=\log(2+R)$. Since $R\ge1$ and $h/\log3 \ge 1$, we have
\begin{align*}
\log(K+2)
&\le
\log(4M_r)
+\frac{2r-1}{r-1}\log R
+\frac{r}{r-1}\log h\\
&\le
\left(
\frac{\log(4M_r)}{\log3}+\frac{3r-1}{r-1}
\right)h.
\end{align*}
Consequently,
\[
(\log(K+2))^{\frac{r}{r-1}}
\le
\left(
\frac{3r-1}{r-1}+\frac{\log(4M_r)}{\log3}
\right)^{\frac{r}{r-1}}
h^{\frac{r}{r-1}}
\le M_rh^{\frac{r}{r-1}},
\]
and therefore,
\[
\expec[\|Tx_K-x_K\|_p]
\le
Q_p\left(
\frac{M_rh^{r/(r-1)}}
{M_rR^{(2r-1)/(r-1)}h^{r/(r-1)}}
\right)^{\frac{r-1}{2r-1}}
=\epsilon.
\]
\end{proof}

In the literature on oracle complexity with iterate-dependent martingale-difference noise,  \cite{tamizholi2026optimal} establishes an $O(\epsilon^{-(4+\eta)})$ oracle complexity for every fixed $\eta>0$ in Euclidean space, and \cite{bravo2024} establishes an $O(\epsilon^{-(10+\eta)})$ oracle complexity for every fixed $\eta>0$ in non-Euclidean spaces using SKM. Our Corollary~\ref{cor_sam} shows that our sample complexity $\widetilde{O}(\epsilon^{-3})$ improves both results for $p \ge 2$, and for $p=\infty$, our sample complexity matches the high-dimensional lower bound of \cite{diakonikolas2026oracle} up to logarithmic factors in $\epsilon$.

\section{Application to stochastic optimization}\label{sec_so}

In this section, we apply dual-anchoring to stochastic gradient descent to efficiently minimize the norm of the gradient.

\subsection{Stochastic optimization}

Consider the expected cost
$F(x)=\mathbb{E}_{z}[f(x,z)]$, where
$f:\mathbb{R}^d\times\cZ\to\mathbb{R}$ and $z$ is a random variable taking values in $\cZ$. The goal of stochastic optimization \cite{nemirovski2009robust} is to solve
\[\min_{x\in \real^d} \mathbb{E}_{z}[f(x,z)].\]
In the case where we can only access $\nabla_x f(x_n,z_{n+1})$ at the given iterate $x_n$ at each iteration $n$, one classical approach to solving the stochastic optimization problem is to apply stochastic gradient descent (SGD) \cite{polyak1992acceleration}:
\begin{equation}
\tag{SGD}
x_n = x_{n-1} - \gamma \nabla_x f(x_{n-1}, z_n),
\end{equation}
 where $z_n$ is an independent copy of $z$ drawn at iteration $n$.
 
 The following assumption is used in the convergence analysis of SGD and corresponds to our Assumption~\ref{assum_fixp}. 
\begin{assumption}\label{assum_sgd}
    $F$ is convex and $L$-smooth with $L>0$. $x_\star\in\argmin_{x\in \real^d} F(x)$. For almost every $z$, $f(\cdot,z)$ is differentiable, and $\nabla F(x)=\mathbb{E}[\nabla_x f(x,z)]$ for every $x$.
\end{assumption}

Under Assumption~\ref{assum_sgd}, if we define
\[
Tx=x-\eta \nabla F(x),
\]
then, by the cocoercivity of $\nabla F$ \cite[Theorem~2.1.5]{nesterov2013introductory},
$T$ is nonexpansive in the Euclidean norm for
$0\leq\eta\leq 2/L$. Then, we have
\[
x_{n}-\eta \nabla_x f(x_{n},z_{n+1})=Tx_{n}+U_{n+1}, \quad \text{where}\quad U_{n+1}=\eta
\left[\nabla F(x_{n})-\nabla_x f(x_{n},z_{n+1})\right],
\]
and $(U_n)_{n\ge1}$ is a martingale difference sequence.

 Now, let $\mathcal F_{k}$ denote the history before $z_{k+1}$ is drawn.
Then $x_{k}$ is $\mathcal F_{k}$-measurable, and $z_{k+1}$ is independent of $\mathcal F_{k}$. 
We consider the following assumption, which is widely used in the analysis of stochastic gradient descent.

\begin{assumption}[ABC condition, \cite{liu2022almost}]\label{assum_abc}
There exist $A,B,C \ge 0$ such that, for every $\mathcal F_k$-measurable $x\in\real^d$,  
\[
\mathbb{E}\!\left[\|\nabla_x f(x,z_{k+1})\|_2^2
\mid\mathcal F_k\right]
\leq A(F(x)-F(x_\star))+B\|\nabla F(x)\|_2^2+C.
\]
\end{assumption}

 Since $\nabla F(x_\star)=0$,
$L$-smoothness gives \cite[Theorem~2.1.5]{nesterov2013introductory}
\[
F(x)- F(x_\star)\leq\frac{L}{2}\|x-x_\star\|_2^2,
\qquad
\|\nabla F(x)\|_2
=\|\nabla F(x)-\nabla F(x_\star)\|_2
\leq L\|x-x_\star\|_2.
\]
Combining these inequalities with the ABC condition and using
$\|x-x_\star\|_2^2\leq 2\|x\|_2^2+2\|x_\star\|_2^2$, we obtain
\[
\mathbb{E}\!\left[\|\nabla_x f(x,z_{k+1})\|_2^2
\mid\mathcal F_k\right]
\leq a+b\|x\|_2^2,\qquad
b=AL+2BL^2,\quad a=C+b\|x_\star\|_2^2.
\]

Furthermore, using
$\mathbb{E}[\nabla_x f(x_{k},z_{k+1})\mid\mathcal F_{k}]=\nabla F(x_{k})$,
we get
\[
\begin{aligned}
\mathbb{E}\!\left[\|U_{k+1}\|_2^2\mid\mathcal F_{k}\right]
&=\eta^2\left(
\mathbb{E}\!\left[\|\nabla_x f(x_{k},z_{k+1})\|_2^2\mid\mathcal F_{k}\right]
-\|\nabla F(x_{k})\|_2^2\right)\\
&\leq\eta^2 a+\eta^2 b\|x_{k}\|_2^2.
\end{aligned}
\]
Therefore, Assumption~\ref{assum_abc} in stochastic optimization reduces to our Assumption~\ref{assum_stoo}.

\subsection{ Stochastic Dual-anchored gradient descent}
Now, we introduce \emph{Stochastic Dual-Anchored Gradient Descent} (SDAGD):
\[
\begin{aligned}
x_{k+\frac12}&=x_k-\alpha\eta\,\nabla_x f(x_k,z_{k+1}),\\
x_{k+1}&=x_k+\frac{K-k}{K-k+2}
              \bigl(x_{k+\frac12}-x_{k-\frac12}\bigr),
\qquad k=0,\ldots,K-1,
\end{aligned}
\]
where $x_{-1/2}=x_0$. By the previous argument, 
\[x_k-\alpha\eta\,\nabla_x f(x_k,z_{k+1}) =(1-\alpha)x_k+\alpha\,(Tx_k+U_{k+1}).\] 
So, by applying Theorem~\ref{thm_resi}, we obtain the following convergence result.
\begin{corollary}
Assume Assumptions~\ref{assum_sgd} and~\ref{assum_abc}. Let $0<\eta\le 2/L$ and $\alpha=D^{-2}(KL_K)^{-2/3}$, where $L_K=2+4\sum_{k=0}^{K-1}(K-k+2)^{-1}$ and  $D=\max\{1,3\eta\sqrt b\}$.
Then SDAGD exhibits the rate
\[
\mathbb E\!\left[\|\nabla F(x_K)\|_2\right]
\le
K^{-\frac13}\bigl(\ln(K+2)\bigr)^{\frac23}M,
\]
where $M=
4^{2/3}\left(\frac{8\kappa D^2}{\eta}+\frac{2}{D}\sqrt{a+b(\|x_\star\|_2+\|x_0-x_\star\|_2)^2}\right)$, with $b=AL+2BL^2$ and $a=C+b\|x_\star\|_2^2$.
For $0<\epsilon\le M$, this implies
\[
\expec[\|\nabla F(x_K)\|_2]
\le\epsilon
\]
with $K=\widetilde O\left((M/\epsilon)^3\right)$ iterations.
\end{corollary}

Note that $\|Tx_K-x_K\|_2=\eta\|\nabla F(x_K)\|_2$. Also, we briefly note that the $\Omega(\epsilon^{-2})$ lower bound in \cite{foster2019complexity} applies to the ABC class with $A,B,C>0$, and \cite{ji2026computation} establishes near-optimal sample complexity under its state-dependent noise condition using accumulative regularization.

\subsection{Simple numerical experiments}
Motivated by a worst-case example in \cite{nesterov2013introductory}, we consider a stochastic regression problem where the  objective function is
\[
\frac12\mathbb E_{a\sim\mathcal N(0,H)}[(a^\top x-y)^2]
=\frac12(x-x_\star)^\top H(x-x_\star),
\]
where   \(H=\tfrac14\operatorname{tridiag}(-1,2,-1)\), $y=a^\top x_\star$ and
$x_\star=\frac{1}{d+1}(d,d-1,\ldots,1)^\top$. Let $F(x) =\mathbb E_{a\sim\mathcal N(0,H)}[f(x,a)]$ and $f(x,a)=\frac12(a^\top x-y)^2$, with \(a,x,x_\star\in\mathbb R^d\).
We consider $d\in\{400,600\}$ and use one fresh sample per iteration.
Since $\nabla F(x)=H(x-x_\star)$ and $0\prec H\preceq I$, $F$ is convex and $1$-smooth.
Thus, $T=I-\eta\nabla F$ is nonexpansive in the Euclidean norm for $0\leq\eta\leq 2$. We also have
\[
\nabla_x f(x,a_k)=a_k\bigl(a_k^\top x-y_k\bigr),
\qquad y_k=a_k^\top x_\star.
\]
For any $\mathcal F_k$-measurable $x$, the stochastic gradient is conditionally unbiased and satisfies
\[
\mathbb E\!\left[\|\nabla_x f(x,a_{k+1})-\nabla F(x)\|_2^2
\mid\mathcal F_k\right]
=
2\operatorname{tr}(H)\bigl(F(x)-F(x_\star)\bigr)
+\|\nabla F(x)\|_2^2.
\]
Consequently, Assumption~\ref{assum_abc} is satisfied with
$A=2\operatorname{tr}(H)$, $B=2$, and $C=0$.
\begin{figure}[h]
    \centering
    \includegraphics[width=1\linewidth]{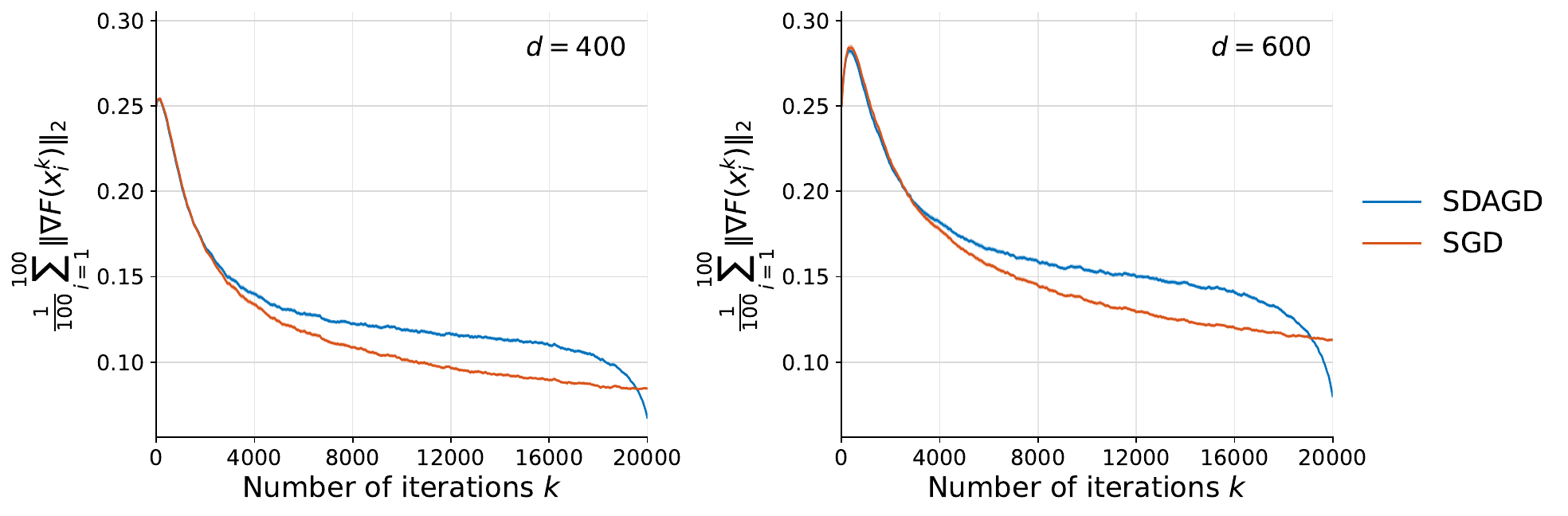}
    \vspace{-.15in}
    \caption{
    Comparison of SDAGD and SGD for $d=400$ (left) and $d=600$ (right).
    Each panel shows the gradient norm across $K=20000$ iterations, averaged over $100$ independent trials.
    At the final iteration, the empirical mean residuals for SDAGD and SGD are
    $0.06745$ and $0.08455$, respectively, for $d=400$,
    and $0.07997$ and $0.11320$, respectively, for $d=600$.
    The shaded regions indicate approximate pointwise $95\%$ confidence intervals for the mean residual.
    }
    \label{fig:_so}
\end{figure}

With the previous objective function, we run SDAGD and SGD for $20000$ iterations for $d=400$ and $d=600$ and record the gradient norm averaged over $100$ independent trials. We set $x_{-1/2}=x_0=0$ for SDAGD and $x_0=0$ for SGD, and use $\alpha=0.0025$  and $\eta=1$. 
Within each trial, both methods use the same regression sample at each iteration. As Figure~\ref{fig:_so} shows,  SGD attains a smaller mean gradient norm during part of the run, but, near the final iteration, SDAGD attains a smaller mean gradient norm in both dimensions.

\section{Application to two-player zero-sum games with the mean payoff}\label{sec_game}

In this section, we apply dual-anchoring to synchronous $Q$-learning to efficiently minimize the span seminorm of the Shapley residual and the duality gap.

\subsection{Two-player zero-sum Markov game}
Consider a two-player zero-sum  Markov game with simultaneous actions, where the max-player tries to maximize the mean payoff, whereas the min-player tries to minimize it over an infinite horizon \cite{mertens1981stochastic}. 

\paragraph{Markov game with the mean payoff}
Let $\Delta(\cX)$ be the space of probability distributions over a set $\cX$. We write $(\cS, \cA^1, \cA^2, P, r)$ to denote the infinite-horizon stochastic game with finite state space $\cS$, finite action spaces $\cA^1, \cA^2$, transition kernel $P\colon \cS \times \cA^1 \times\cA^2 \rightarrow \Delta(\cS)$, and bounded reward $r\colon \cS \times \cA^1\times \cA^2 \rightarrow [-B,B]$. Let $\pi^i\colon \cS \rightarrow \Delta(\cA^i)$ denote a stationary policy for player $i$, $i=1,2$, and define
\begin{align*}
    J^{\pi^1,\pi^2}(s)=\lim_{T\rightarrow \infty} \frac{1}{T}\expec_{\pi^1,\pi^2}\bigg[\sum^{T-1}_{t=0} r(s_t,a_t^1,a_t^2) \,\Big|\, s_0=s\bigg]
\end{align*}
as the mean payoff of a given policy pair, where $\expec_{\pi^1,\pi^2}$ denotes the expected value over all trajectories $(s_0, a^1_0, a^2_0, s_1, a^1_1,a^2_1, \dots, s_{T-1}, a^1_{T-1}, a^2_{T-1})$ induced by $P$ and $\pi^1,\pi^2$, with the players' actions sampled independently conditional on the current state. The limit exists for every stationary policy pair because the Ces\`aro limit of every finite stochastic matrix exists \cite[Theorem A.6]{Puterman2014}. Write $r_{\pi^1,\pi^2}(s)=\mathbb{E}_{a \sim \pi^1(\cdot\,|\,s), b \sim \pi^2(\cdot\,|\,s)}\left[r(s,a,b)\right]$ for the reward induced by the policy pair $\pi^1,\pi^2$. 
For simplicity, we write $\cA=\cA^1\times\cA^2$ and $d=|\cS||\cA^1||\cA^2|$.

Next, define the Shapley operator \cite{shapley1953stochastic} by
\[
TQ(s,a,b)
=
r(s,a,b)+\sum_{s' \in \cS}P(s'|s,a,b)\max_{\mu\in\Delta(\cA^1)}
\min_{\nu\in\Delta(\cA^2)}
\sum_{\substack{a' \in \cA^1 \\ b' \in \cA^2}}\mu(a')\nu(b')Q(s',a',b'),
\]
for $Q\in\real^d$. With this Shapley operator $T$, we consider the following assumption.

\begin{assumption}[Ergodic equation]\label{assum_erg}
There exist $Q^\star\in\real^{|\cS||\cA|}$ and $J^\star\in\real$ such that $TQ^\star=Q^\star+J^\star\mathbf1$.
\end{assumption}
In the zero-sum game literature,
this condition has been studied in connection with the existence of solutions to zero-sum stochastic games \cite{akian2014ergodicity} and also is more general than the communication condition in \cite[Assumption 2 and Theorem E.1]{wei2017online}.
 Under Assumption~\ref{assum_erg}, 
the existence of optimal stationary policies is guaranteed and $J^\star$ is equal to the optimal mean payoff of the zero-sum game (\cite[Proposition 2.1]{akian2014ergodicity} and \cite[Lemma 1]{federgruen1980successive}). Thus,  
under  Assumption~\ref{assum_erg}, the optimal mean payoff satisfies
\[
J^\star = \max_{\pi^1}\min_{\pi^2}
J^{\pi^1,\pi^2}(s),\qquad \forall s\in\cS.
\]

Next, to formulate the two-player zero-sum game as a fixed-point equation, we consider the following  quotient space. Define the quotient space $\hat{\real}^{\cS\times\cA}=\real^{\cS\times\cA}/E$, where $E=\{c\mathbf1:c\in\real\}$ is the subspace of constant vectors, and endow the quotient space with
the induced quotient norm
$$\big\|[Q]\big\|_E=\min_{c\in\real}\big\|Q+c\mathbf1\big\|_\infty.$$
Then, the Shapley operator $T:\real^d\to\real^d$ induces the operator
$\widehat T([Q])=[TQ]$ on the quotient space, since $T(Q+c\mathbf1)=TQ+c\mathbf1$. Also, one can readily check that  $[Q]\in\operatorname{Fix}(\widehat T)$ if and only if $TQ-Q\in E$.  
Then, Assumption~\ref{assum_erg} gives $\widehat T([Q^\star])=[Q^\star]$ on the quotient space. 

Now, define the span seminorm $\sp{\cdot}$ by $\sp{x}=\max_j x_j-\min_j x_j$ for a vector $x$. Then we get
 \[\big\|[Q]\big\|_E=\frac12\sp{Q},\]
and one can check that the Shapley operator $T$ is nonexpansive with respect to $\sp{\cdot}$, and this also implies $\widehat T$ is nonexpansive with respect to $\|\cdot\|_E$.


Overall, in this quotient space, the ergodic equation gives the fixed-point equation $\widehat T([Q^\star])=[Q^\star]$ since $J^\star \mathbf{1}$ is a constant vector. This corresponds to our Assumption~\ref{assum_fixp}, and the quotient norm can be calculated from the span seminorm.

Based on these, now, we consider the synchronous sampling setup. Instead of accessing the Shapley operator $T$ directly, for a given $Q$, suppose we can access \[
(\widetilde TQ)(s,a,b)
=
r(s,a,b)+\max_{\mu\in\Delta(\cA^1)}
\min_{\nu\in\Delta(\cA^2)}
\sum_{a' \in \cA^1,b'\in \cA^2}\mu(a')\nu(b')Q(s',a',b')
,
\]
for all $s,a,b$, with $s' \sim P(\cdot|s,a,b)$.  One classical approach to solving the fixed-point equation in this setup is synchronous Q-learning \cite{littman1994markov}:
\[
\begin{multlined}
Q_{k+1}(s,a,b)=(1-\alpha)Q_k(s,a,b)\\
{}+\alpha\left(r(s,a,b)+\max_{\mu\in\Delta(\cA^1)}
\min_{\nu\in\Delta(\cA^2)}
\sum_{a',b'}\mu(a')\nu(b')Q_k(s',a',b')\right)\quad\forall(s,a,b),
\end{multlined}
\]
which is equivalent to 
\[Q_{k+1}(s,a,b)=(1-\alpha)Q_k(s,a,b) + \alpha(TQ_k+U_{k+1})(s,a,b),\]
where
\[U_{k+1}(s,a,b)=r(s,a,b)+\max_{\mu\in\Delta(\cA^1)}
\min_{\nu\in\Delta(\cA^2)}
\sum_{a',b'}\mu(a')\nu(b')Q_k(s',a',b')-(TQ_k)(s,a,b).\]
Then we can directly check that $(U_k)_{k\ge1}$ is a martingale difference sequence. Moreover, the following fact for the span seminorm follows from \cite[Section 2.2, Example 3.2, and Proposition 3.4]{juditsky2008large} using the previously explained equivalence between the span seminorm and the quotient norm.

\begin{fact}[\cite{juditsky2008large}]\label{fact_sam_sp}
Let
$U_1,\ldots,U_m\in\real^d$ be a martingale difference sequence
with respect to a filtration $(\mathcal F_i)_{i=0}^m$, satisfying
\[
\expec[U_i\mid\mathcal F_{i-1}]=0,
\qquad
\expec[\sp{U_i}^2]<\infty,
\qquad i=1,\ldots,m.
\]
Then, for any deterministic coefficients $a_1,\ldots,a_m\in\real$,
\[
\left(
\expec\left[
\sp{\sum_{i=1}^m a_i U_i}^{2}
\right]
\right)^{1/2}
\le
C
\left(
\sum_{i=1}^m
|a_i|^{2}
\expec\left[\sp{U_i}^{2}\right]
\right)^{1/2},
\]
where $C=\sqrt{2e\ln(2d)}$.
\end{fact}

Since $Q_k$ is $\cF_k$-measurable and the samples have the specified conditional distributions, we have $\expec[U_{k+1}\mid\cF_k]=0$, and by direct calculations, 
\begin{align*}
\sp{U_{k+1}}^2&\le4\sp{Q_k}^2\le(\sp{r}+2\sp{Q_k})^2\le2\sp{r}^2+8\sp{Q_k}^2.
\end{align*}
Thus, Assumption~\ref{assum_stoo} holds in the span seminorm with $a=2\sp{r}^2$ and $b=8$.

\subsection{Dual-Anchored Q-learning}

Now we are ready to present our \emph{Dual-Anchored Q-learning} (DA-$Q$-learning) with synchronous updates:
\begin{align*}
Q_{k+\frac{1}{2}}(s,a,b)&=(1-\alpha)Q_k(s,a,b)\\&+\alpha\biggl(r(s,a,b)+\max_{\mu\in\Delta(\cA^1)}\min_{\nu\in\Delta(\cA^2)}
\sum_{a',b'}\mu(a')\nu(b')Q_k(s',a',b')\biggr)\qquad \forall(s,a,b),\\
Q_{k+1}
&=Q_k+\frac{K-k}{K-k+2}(Q_{k+\frac{1}{2}}-Q_{k-\frac{1}{2}})
\qquad 0\le k\le K-1,
\end{align*}
where, conditional on the history $\cF_k$, $s' \sim P(\cdot\mid s,a,b)$ for each $(s,a,b)$ at iteration $k$.
$Q_{-\frac12}=Q_0$, and the constant $\alpha$ is chosen as in the corollary below.

 By the definition of $U_{k+1}$, we get
\begin{align*}
r(s,a,b)+\max_{\mu\in\Delta(\cA^1)}
\min_{\nu\in\Delta(\cA^2)}
\sum_{a',b'}\mu(a')\nu(b')Q_k(s',a',b')=(TQ_k+U_{k+1})(s,a,b).
\end{align*}
Hence, we can apply the proof of Theorem~\ref{thm_resi} on the quotient space and, using the equivalence between the quotient norm and the span seminorm, we can use Fact~\ref{fact_sam_sp} to obtain the sample complexity.

Before presenting the convergence result for DA-$Q$-learning, we introduce the duality gap \cite{nemirovski2004prox,aubin2007mathematical}, a measure widely used in zero-sum games. Define the duality gap as
\[
\operatorname{Gap}(\pi^1,\pi^2)
=
\|\max_{\pi} J^{\pi,\pi^2}
-
\min_\pi J^{\pi^1,\pi}\|_\infty.
\]
Note that the duality gap bounds both players' optimality errors as follows:
\[
\left\|\sup_{\pi}J^{\pi,\pi^2}-J^\star\mathbf1\right\|_\infty
\le\operatorname{Gap}(\pi^1,\pi^2),\qquad
\left\|J^\star\mathbf1-\inf_{\pi}J^{\pi^1,\pi}\right\|_\infty
\le\operatorname{Gap}(\pi^1,\pi^2).
\]
Also, the duality gap can be bounded by the span seminorm of the Shapley residual as the following lemma shows. We defer the proof to the Appendix.
\begin{lemma}\label{erg_res} Under Assumption~\ref{assum_erg},
    \[\operatorname{Gap}(\pi^1_Q,\pi^2_Q)\le\sp{TQ-Q},\]
    where $(\pi^1_Q(s),\pi^2_Q(s))$ is induced by $TQ$ at each state $s$.
\end{lemma}

Now, we are ready to present the convergence result for DA-$Q$-learning.

\begin{corollary}
Assume Assumption~\ref{assum_erg}. Let $
\alpha=\frac{(KL_K)^{-\frac23}}{144e\ln(2d)}$,
 where  $d=|\cS||\cA^1||\cA^2|$ and
$L_K=2+4\sum_{k=0}^{K-1}(K-k+2)^{-1}$. Dual-Anchored $Q$-learning exhibits the rate
\[
\expec[\operatorname{Gap}(\pi^1_{Q_K},\pi^2_{Q_K})]
\le\mathbb E\!\left[\sp{TQ_K-Q_K}\right]
\le
K^{-\frac13}\bigl(\ln(K+2)\bigr)^{\frac23}M,
\]
where $M=
8\times10^3\ln(2d)\sp{Q_0-Q^\star}
+\sqrt{
\sp{r}^2
+4\left(\sp{Q^\star}+\sp{Q_0-Q^\star}\right)^2
}$.
For $0<\epsilon\le M$, this implies
\[
\expec[\operatorname{Gap}(\pi^1_{Q_K},\pi^2_{Q_K})]
\le
\expec[\sp{TQ_K-Q_K}]
\le\epsilon
\]
with $K=\widetilde O\left((M/\epsilon)^3\right)$ iterations,
requiring $dK$ transition samples.
\end{corollary}


While prior work has studied regret for model-based methods \cite{wei2017online} and asymptotic convergence of
$Q$-learning methods \cite{mannor2004reinforcement} for mean-payoff zero-sum Markov games,
to the best of our knowledge, our work provides the first
explicit last-iterate sample-complexity analysis in a synchronous sampling setup under the ergodic equation assumption.

\begin{remark}
We note that our results also apply to single-agent average-reward reinforcement learning. In the single-agent case, every finite weakly communicating MDP satisfies Assumption~\ref{assum_erg} \cite[Section 9]{Puterman2014}. So, when $\cA^2=\{a\}$, the expected duality gap reduces to the expected policy error $
\expec[\operatorname{Gap}(\pi^1_{Q_K},\pi^2_{Q_K})]
=
\expec\left[\left\|J^\star\mathbf1-J^{\pi^1_{Q_K}}\right\|_\infty\right].$ Therefore, DA-$Q$-learning also gives $\widetilde{O}(\epsilon^{-3})$ iterations to obtain an expected policy error of at most $\epsilon$ for average-reward reinforcement learning.
\end{remark}

\subsection{Simple numerical illustration}
Motivated by the worst-case MDP construction in Theorem~3 of
\cite{lee2025multi} for average-reward reinforcement learning, we consider
the following zero-sum Markov game $(\cS,\cA_1,\cA_2,P,r)$ with one action per player:
\begin{align*}
    &\cS=\{s_1,\dots,s_n\}, \qquad
    \cA_1=\{a_1\}, \qquad
    \cA_2=\{a_2\}, \qquad
    r(s_i,a_1,a_2)=\mathbf{1}_{\{i=1\}},\\
    &P(s_i\mid s_j,a_1,a_2)
    =(1-\delta)
    \mathbf{1}_{\{(i,j)=(n-1,1)\ \text{or}\ j=i+1\}}
    +\delta\mathbf{1}_{\{i=j\}},
\end{align*}
where $n=1002$ and $\delta=0.1$.
We set $x_0=0$ and perform $K=1000$ synchronous iterations,
sampling one successor independently for each state at each iteration.

\begin{figure}[h]
    \centering
    \includegraphics[width=1\linewidth]{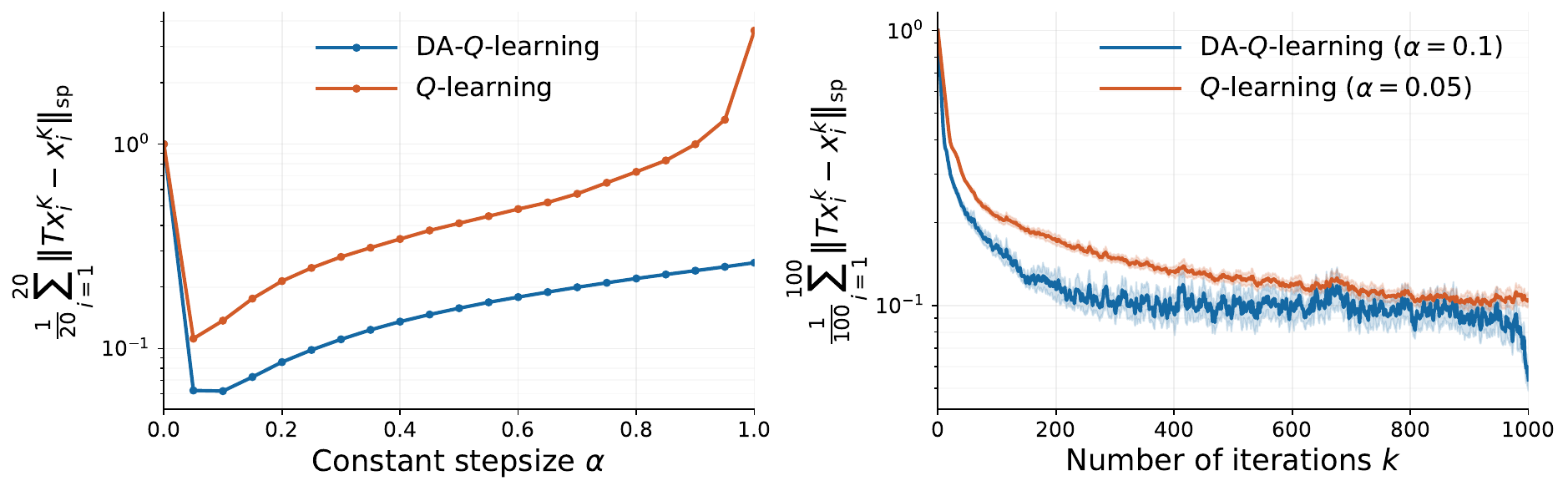}
        \vspace{-.15in}
    \caption{
    Comparison of DA-$Q$-learning and $Q$-learning.
    Left: final empirical span residuals for each tested constant stepsize
    $\alpha\in[0,1]$, averaged over $20$ independent trials.
    Right: empirical span residuals across iterations, averaged over
    $100$ independent trials, using $\alpha=0.1$ for DA-$Q$-learning and
    $\alpha=0.05$ for  $Q$-learning.
    At $K=1000$, the empirical residuals
    $\frac{1}{100}\sum_{i=1}^{100}
    \|T x_i^K-x_i^K\|_{\mathrm{sp}}$
    are approximately $0.05372$ for DA-$Q$-learning and $0.10448$ for  $Q$-learning.
    The shaded regions in the right panel indicate approximate pointwise
    $95\%$ confidence intervals for the mean residual.
    }
    \label{fig_zs}
\end{figure}

To examine the effect of dual-anchoring, we compare $Q$-learning and DA-$Q$-learning
with constant stepsizes $\alpha$, keeping $m_k$ as in the theorem.
First, we evaluate the final empirical
span residual on a uniform grid of stepsizes in $[0,1]$,
averaging over $20$ independent trials for each stepsize.
As the left panel in Figure~\ref{fig_zs} shows, DA-$Q$-learning attained a smaller empirical mean residual than $Q$-learning.
Using the selected stepsizes $\alpha=0.1$ for DA-$Q$-learning and $\alpha=0.05$
for $Q$-learning, we then record the empirical span residual across iterations
over $100$ independent test trials. As shown in the right panel of Figure~\ref{fig_zs},
DA-$Q$-learning achieves approximately half the residual
of $Q$-learning at the last iterate.

\section{Conclusion}
In this work, we showed that SDAKM achieves efficient oracle complexity for nonexpansive operators under iterate-dependent martingale-difference noise in $(\mathbb{R}^d,\|\cdot\|_p)$. In particular, SDAKM utilizes the averaging effect in SKM, which efficiently reduces stochastic errors, and dual-anchoring in Dual Halpern iteration, which reduces the term involving the distance to the solution. Based on these mechanisms, SDAKM achieves near-optimal oracle complexity for $p=\infty$, up to logarithmic factors in $\epsilon$, and improves the state-of-the-art sample complexity in both Euclidean and non-Euclidean spaces.


One interesting direction for future research is developing an anytime variant of SDAKM that achieves almost sure convergence as $K\rightarrow\infty$. Establishing lower bounds for $p \in (2,\infty)$ to check the optimality of our results would also be interesting. Lastly, for two-player zero-sum games, extending the method to asynchronous settings with additional analysis would also be a promising research direction.


\paragraph{Disclosure of AI Use}
In this work, the author used AI tools for refining calculations and implementing the numerical experiments. 
We have not used generative AI tools for the development of the core methodology and  theoretical results. The author assumes responsibility for all content.

\newpage

\bibliographystyle{plain}
\bibliography{main}

\newpage

\section{Appendix}\label{app}
We present the omitted proof of Lemma~\ref{erg_res}, which adapts \cite[Proposition 2.1]{lee2025near}.
\begin{proof}
For any stationary pair $(\pi^1,\pi^2)$, let $P_{\pi^1,\pi^2}$ be its induced transition matrix, and define its Ces\`aro limit by
\[
P^\star_{\pi^1,\pi^2}
:=
\lim_{N\to\infty}\frac1N\sum_{t=0}^{N-1}P_{\pi^1,\pi^2}^{\,t}.
\]
Since the state space is finite, this limit exists and is a stochastic matrix satisfying $P^\star_{\pi^1,\pi^2}P_{\pi^1,\pi^2}=P^\star_{\pi^1,\pi^2}$.
Then, $J^{\pi^1,\pi^2}=P^\star_{\pi^1,\pi^2}r_{\pi^1,\pi^2}.$
Also, for an arbitrary $Q$, set
\[
\begin{gathered}
V_Q(s)=\max_{\pi^1}
\min_{\pi^2}
\!\!\!\!\!\!\sum_{a\in \cA^1,b\in \cA^2}\!\!\!\!\!\!\pi^1(a|s)\pi^2(b|s)Q(s,a,b),\\ e=TQ-Q,\quad
\ell=\min_{s,a,b}e(s,a,b),\quad
u=\max_{s,a,b}e(s,a,b),
\end{gathered}
\]
and let $\pi_Q^1$ and $\pi_Q^2$ be the policies induced by $TQ$ at each state $s$.

For every stationary policy $\pi$ of the min-player, $Q_{\pi_Q^1,\pi}\ge V_Q$ holds,
where $Q_{\pi_Q^1,\pi}(s)$ denotes the expectation of $Q(s,a,b)$ under the two policies. Define $e_{\pi_Q^1,\pi}(s)$ analogously. Then, we have
\[
r_{\pi_Q^1,\pi}+P_{\pi_Q^1,\pi}V_Q-V_Q
=Q_{\pi_Q^1,\pi}-V_Q+e_{\pi_Q^1,\pi}
\ge\ell\mathbf1.
\]
Multiplying by $P^\star_{\pi_Q^1,\pi}$ preserves this inequality, and by $P^\star_{\pi^1,\pi^2}P_{\pi^1,\pi^2}=P^\star_{\pi^1,\pi^2}$, we get
\[
P^\star_{\pi_Q^1,\pi}r_{\pi_Q^1,\pi}
+\bigl(P^\star_{\pi_Q^1,\pi}P_{\pi_Q^1,\pi}
-P^\star_{\pi_Q^1,\pi}\bigr)V_Q\ge
\ell\,P^\star_{\pi_Q^1,\pi}\mathbf1.
\]
Consequently,
$J^{\pi_Q^1,\pi}\ge\ell\mathbf1$ for every stationary policy $\pi$ of the min-player.

Similarly, for every stationary policy $\pi$ of the max-player,  $Q_{\pi,\pi_Q^2}\le V_Q$, and we get
\[
r_{\pi,\pi_Q^2}+P_{\pi,\pi_Q^2}V_Q-V_Q
=Q_{\pi,\pi_Q^2}-V_Q+e_{\pi,\pi_Q^2}\le u\mathbf1.
\]
Multiplying by  $P^\star_{\pi,\pi_Q^2}$ yields
$J^{\pi,\pi_Q^2}\le u\mathbf1$ for every stationary policy $\pi$ of the max-player.

Thus, for every initial state $s$,
\[
\begin{aligned}
\sup_\pi J^{\pi,\pi_Q^2}(s)
-\inf_\pi J^{\pi_Q^1,\pi}(s)&\le u-\ell=\sp{TQ-Q}.
\end{aligned}
\]
Taking the maximum over $s$ proves the claim.
\end{proof}

\end{document}